\documentclass[11pt]{amsart}

\usepackage{amsmath,amsfonts,amsthm,amsopn,amssymb}
\usepackage{verbatim,wasysym,mathrsfs}
\usepackage[left=2.6cm,right=2.6cm,top=3.1cm,bottom=3.1cm]{geometry}
\usepackage{cite,marginnote}
\usepackage{color,enumitem,graphicx}
\usepackage[colorlinks=true,urlcolor=blue, citecolor=red,linkcolor=blue,
linktocpage,pdfpagelabels, bookmarksnumbered,bookmarksopen]{hyperref}
\usepackage[english]{babel}
\usepackage[T1]{fontenc}
\usepackage{lmodern}
\usepackage[hyperpageref]{backref}

\numberwithin{equation}{section}

\makeindex
\newtheorem{theorem}{Theorem}[section]
\newtheorem{proposition}[theorem]{Proposition}
\newtheorem{lemma}[theorem]{Lemma}
\newtheorem{remark}[theorem]{Remark}
\newtheorem{example}[theorem]{Example}
\newtheorem{corollary}[theorem]{Corollary}
\newtheorem{definition}[theorem]{Definition}

\newcommand{\ud}{\mathrm{d}}
\newcommand{\RN}{\mathbb R^N}

\newcommand{\iy}{\infty}

\newcommand{\s}{\section}

\newcommand{\DD}{\Delta}

\newcommand{\na}{\nabla}

\newcommand{\la}{\lambda}

\newcommand{\pa}{\partial}

\newcommand{\R}{\mathbb R}

\newcommand{\ti}{\tilde}

\newcommand{\re}[1]{(\ref{#1})}

\newcommand{\rg}{\rightarrow}

\newcommand{\lan}{\langle}
\newcommand{\ran}{\rangle}
\newcommand{\e}{\varepsilon}
\newcommand{\vp}{\varphi}

\newcommand{\lab}{\label}
\newcommand{\bt}{\begin{theorem}}
\newcommand{\et}{\end{theorem}}
\newcommand{\bl}{\begin{lemma}}
\newcommand{\el}{\end{lemma}}
\newcommand{\bd}{\begin{definition}}
\newcommand{\ed}{\end{definition}}
\newcommand{\bc}{\begin{corollary}}
\newcommand{\ec}{\end{corollary}}
\newcommand{\bp}{\begin{proof}}
\newcommand{\ep}{\end{proof}}
\newcommand{\bx}{\begin{example}}
\newcommand{\ex}{\end{example}}
\newcommand{\bi}{\begin{exercise}}
\newcommand{\ei}{\end{exercise}}
\newcommand{\bo}{\begin{proposition}}
\newcommand{\eo}{\end{proposition}}
\newcommand{\br}{\begin{remark}}
\newcommand{\er}{\end{remark}}
\newcommand{\be}{\begin{equation}}
\newcommand{\ee}{\end{equation}}
\newcommand{\ba}{\begin{align}}
\newcommand{\ea}{\end{align}}
\newcommand{\bn}{\begin{enumerate}}
\newcommand{\en}{\end{enumerate}}
\newcommand{\bg}{\begin{align*}}
\newcommand{\bcs}{\begin{cases}}
\newcommand{\ecs}{\end{cases}}

\def\R{\mathbb R}

\makeatletter
\def\@makefnmark{}
\makeatother

\newcommand{\pl}{\partial}

\newcommand{\bean}{\begin{eqnarray*}}
\newcommand{\eean}{\end{eqnarray*}}

\renewcommand{\epsilon}{\varepsilon}

\title[Modulational stability of ground states]{Modulational stability of ground states \\ to nonlinear Kirchhoff equations}

\author[J. J. Zhang]{Jianjun Zhang}
\author[Z. S. Liu]{Zhisu Liu}
\author[M. Squassina]{Marco Squassina}

\address[J. J.\ Zhang]{\newline\indent College of Mathematics and Statistics,
Chongqing Jiaotong University
\newline\indent
Chongqing 400074, PR China}
\email{\href{mailto:zhangjianjun09@tsinghua.org.cn}{zhangjianjun09@tsinghua.org.cn}}

\address[Z. S. Liu]{\newline\indent School of Mathematics and Physics,
University of South China
\newline\indent
Hengyang, Hunan 421001, P.R. China}
\email{\href{mailto:liuzhisu183@sina.com}{liuzhisu183@sina.com}}

\address[M.\ Squassina]{\newline\indent Dipartimento di Matematica e Fisica,
 Universita Cattolica del Sacro Cuore,
\newline\indent
Via Musei 41, 25121 Brescia, Italy}
\email{\href{mailto:marco.squassina@dmf.unicatt.it}{marco.squassina@unicatt.it}}

\thanks{The corresponding author is M. Squassina.\\ Z. Liu is supported by the NSFC (11626127).
M.\ Squassina is member of the Gruppo Nazionale
	per l'Analisi Matematica, la Probabilit\`a
	e le loro Applicazioni (GNAMPA) of the Istituto Nazionale di Alta Matematica (INdAM).
J. J.\ Zhang was partially supported by NSFC(No.11871123).}


\begin{document}

\begin{abstract}
We investigate the stability of ground states to a nonlinear focusing Schr\"odinger equation in presence of a Kirchhoff term. Through a spectral analysis of the linearized operator about ground states, we show a modulation stability estimate of ground states in the spirit of one due to Weinstein [{\it SIAM J. Math. Anal.}, 16(1985),472-491].
\end{abstract}
\maketitle

\section{Introduction and main result}\label{sec1}

\subsection{Overview}
Let us consider the following nonlinear focusing Kirchhoff equation with a potential and an initial datum
\begin{equation}\lab{q1}
\left\{
\begin{aligned}
    &i\e\pl_tu^\e=-\frac{1}{2}\left(\e^2+\e\int_{\R^3}|\na u^\e|^2\right)\DD u^\e+V(x)u^\e-|u^\e|^{2p}u^\e, & \quad t>0,\,\,\,x\in\R^3, \\
    &u^\e(x,0)=u_0(x), & \quad x\in\R^3,
  \end{aligned}
\right.
\end{equation}
where $p\in(0,2)$ and $\e>0$(referring to Plank's constant). Similar to \cite[Theorem 6.1.1 and Corollary 6.1.2]{Caz}, problem \re{q1} is globally well-posed, provided that $V\in L^{m}(\R^3)+L^\iy(\R^3)$ for some $m>3/2$. Here we refer to \cite{Kirchhoff83} for the background of Kirchhoff equations. Of particular interest are 
the {\it standing wave solutions} of \re{q1}, namely special solutions of \re{q1} of the form 
$$
u^\e(x,t)=v_\e(x)e^{\frac{i}{\e}\theta t},\qquad x\in\R^3,\,\, t\in\R^+,\,\,\theta\in\R.
$$
In this case, $v_\e$ is a solution of the following singularly perturbed Kirchhoff equation
\be\lab{ks}
-\frac{1}{2}\left(\e^2+\e\int_{\R^3}|\na v_\e|^2\right)\DD v_\e+\ti{V}(x)v_\e=|v_\e|^{2p}v_\e, \,\,\,\,x\in\R^3,
\ee 
where $\ti{V}(x)=V(x)+\theta$. An interesting class of solutions to \re{ks} are families of solutions which develop a spike shape around certain points (such as local minima, local maxima and degenerate or non-degenerate critical points of $\ti{V}$) in $\R^3$ as $\e\rg 0$.
These standing wave solutions are very often referred as the semiclassical states for $\e$ small. Initiated by Floer and Weinstein \cite{F-W} for the Schr\"odinger equations $$-\e^2\DD v+W(x)v=f(v),$$ semiclassical states have attracted a considerable attention in the last three decades. For the progress on this topic, we refer e.g.\ to Ambrosetti and Malchiodi \cite{Am} and the reference therein. In the study of the perturbed problem \re{ks}, the following limit problem plays a crucial role
\be\lab{lim}
-\frac{1}{2}\left(1+\int_{\R^3}|\na u|^2\right)\DD u+u=|u|^{2p}u,\quad \,u\in H^1(\R^3).
\ee
It is shown in \cite{Li} that the positive solution of \re{lim} is, up to translation, unique. Denote by $r$ the positive, radially symmetric solutions of \re{lim}. 
\vskip0.1in
Another related topic is to associate to \re{q1} a family of initial data $u_0$ which oscillate or concentrate with scale $\e$, and investigate the evolution of $u^\e$ in time. Precisely, by choosing a suitable initial datum $u_0$ related to the ground state solution $r$, it can be expected that the evolution $u^\e$ remains close to $r$ locally uniformly in time in the semiclassical regime of $\e$ going to zero, driven around by a Newtonian law associated with the potential $V$. This kind of asymptotic behavior is called in the literature {\it soliton dynamics}. In this aspect, we refer the readers to a survey \cite{Tao}. In \cite{Bron}, Bronski and Jerrard considered the following focusing Schr\"odinger equation with a potential
\begin{equation}\lab{bro}
i\e\pl_tu^\e=-\frac{\e^2}{2}\DD u^\e+V(x)u^\e-|u^\e|^{2p}u^\e,  \quad t>0,\,\,\,x\in\RN.
\end{equation}
By using the conservation law (quantum and classical) and the stability of the ground state $Q$ to the limit problem 
$$
-\frac{1}{2}\DD u+u=|u|^{2p}u,\quad u\in H^1(\RN),
$$ 
they proved the solution of \re{bro} exhibits the asymptotic soliton dynamics if the initial datum has the form of 
$$
Q\left(\frac{x-x_0}{\e}\right)e^{i\frac{x\cdot \xi_0}{\e}},\quad x_0, \xi_0\in\RN.
$$
In other words the solution behaves as
$$
Q\left(\frac{x-x(t)}{\e}\right)e^{i\frac{x\cdot \xi(t)}{\e}}
$$
where the parameters $(x(t),\xi(t))$ satisfy the Newton type equation
$$
\frac{\ud x}{\ud t}=\xi(t),\,\,\,\,\,\frac{\ud \xi}{\ud t}=-\na V(x(t)),\quad x(0)=x_0,\,\,\xi(0)=\xi_0.
$$
Subsequently, in \cite{Kera}, Keraani refined and improved the method introduced by Bronski and Jerrard \cite{Bron}.
Later, in \cite{Sel}, Selvitella turned to study the Schr\"odinger equations 
$$
i\pl_tu^\e=-\frac{1}{2}\left(\frac{\e}{i}\na-A(x)\right)^2u^\e+V(x)u^\e-|u^\e|^{2p}u^\e, \quad t>0,\,\,\,x\in\RN
$$
with external electric and magnetic field $V$ and $A$, $B=\na\times A$. Using the linearization argument, the author adopted the idea due to Bronski and Jerrard \cite{Bron} to show the asymptotic evolution of the semiclassical limit as $\e\rg0$. In \cite{Marco1}, Squassina extended and improved the results in \cite{Sel}. 

For more progress in this direction, we also would like to cite \cite{Marco} for nonlocal Choquard equations and \cite{Mon} for systems of weakly coupled Schr\"odinger equations.

\subsection{Main result} In the works above, the nonlinear term is subcritical, namely, $p<2/N$, where $N$ is the dimension. It is well known that the ground states of the associated limit problems above are orbitally stable when $p<2/N$. For more details, we refer the readers to \cite{Lions,Caz}. In the present paper, we also consider the subcritical case: $0<p<2/3$. Moreover, we should point out that in the works above, to establish the soliton dynamics of semiclassical states on finite time intervals, some kind of energy convexity plays an important role. More precisely, via a delicate spectral analysis of the linearized operator at the ground state of the limit problem \re{lim}, we establish a modulational stability result in term of Kirchhoff problems \re{q1}. 

Our main result can read as follows.
\bt\lab{Th1}
There exists $C>0$ such that for any $\phi\in H^1(\R^3,\mathbb{C})$, there holds that
$$
\mathcal{E}(\phi)-\mathcal{E}(r)\ge C\inf_{(x,\theta)\in\R^3\times[0,2\pi)}\|\phi-e^{i\theta}r(\cdot-x)\|^2+o\left(\inf_{(x,\theta)\in\R^3\times[0,2\pi)}\|\phi-e^{i\theta}r(\cdot-x)\|^2\right),
$$
provided that $\|\phi\|_2=\|r\|_2$ and
$$
\inf_{(x,\theta)\in\R^3\times[0,2\pi)}\|\phi-e^{i\theta}r(\cdot-x)\|\le\|r\|.
$$
\et

With the help of Theorem \ref{Th1}, the evolution $u^\e$ of \re{q1} should remain close to $r$ locally uniformly in time, provided a suitable initial datum $u_0$ related to the ground state solution $r$. We will subsequently deal with this topic for the Kirchhoff problem \re{q1}.

\vskip0.1in
\noindent{\bf Notations.}

\begin{itemize}
\item [$\bullet$] For any $z\in\mathbb{C}$, $\bar{z}, \mathfrak{Re}(z), \mathfrak{Im}(z)$ denote the complex conjugate, real and imaginary part.
\item [$\bullet$] For any $z,w\in\mathbb{C}$, it holds that $\mathfrak{Re}(\bar{z}w)=\mathfrak{Re}(z\bar{w})$ and $\mathfrak{Im}(\bar{z}w)=-\mathfrak{Im}(z\bar{w})$.
\item [$\bullet$] For any $z,w\in\mathbb{C}$, we define $z\cdot w=\mathfrak{Re}(z\bar{w})=\frac{1}{2}(z\bar{w}+\bar{z}w)$.
\item [$\bullet$] For any $x,y\in\R^3$, we denote by $x\cdot y$ the inner product between $x$ and $y$.
\item [$\bullet$] $c,C$ denote (possibly different) positive constants which may change from line to line.
\item [$\bullet$] $H^1(\R^3)=H^1(\R^3,\R)$ and $H^1(\R^3,\mathbb{C})$ are real and complex Hilbert space respectively, endowed with the norm
$$
\|u\|=\left(\frac{1}{2}\|\na u\|_2^2+\|u\|_2^2\right)^{\frac{1}{2}},\,\,u\in H^1(\R^3,\mathbb{C}).
$$
\item [$\bullet$] Denote by $(u,v)$ the scalar product in $L^2(\R^3,\mathbb{C})$ and
$$
(u,v)_{H^1}=(u,v)+\frac{1}{2}(\na u,\na v),\,\,\mbox{for}\,\,\, u,v\in H^1(\R^3,\mathbb{C}).
$$
\end{itemize}

\s{Preliminary results}

In this section, we give a few basic properties about the ground state solutions to problem \re{lim}.
\subsection{The limit problem}
It is shown in \cite{Li} that $r$ is the unique radially symmetric solution of \re{lim}. Moreover, it is non-degenerate in the sense that
$$
{\rm Ker}L_+={\rm span}\left\{\pl_{x_1}r,\pl_{x_2}r,\pl_{x_3}r\right\},
$$
where $L_+$ is given as follows
$$
L_+\vp=-\frac{1}{2}\left(1+\int_{\R^3}|\na r|^2\right)\DD\vp-\left(\int_{\R^3}\na r\na\vp\right)\DD r+\vp-(2p+1)r^{2p}\vp,\,\,\vp\in L^2(\R^3).
$$
Moreover, $r\in C^\iy(\R^3)$, $r(0)=\max_{x\in\R^3}r(x)$ and $r,|\na r|$ exponentially decay at infinity.

Now, we consider the following minimization problem with a constraint.
For any $u\in H^1(\R^3,\mathbb{C})$, let
$$
\mathcal{E}(u)=\frac{1}{2}\int_{\R^3}|\na u|^2+\frac{1}{4}\left(\int_{\R^3}|\na u|^2\right)^2-\frac{1}{p+1}\int_{\R^3}|u|^{2p+2},
$$
and
\be\lab{mini}
e:=\inf_{u\in\mathcal{M}}\mathcal{E}(u),\,\,\,\mathcal{M}:=\left\{u\in H^1(\R^3,\mathbb{C}):\|u\|_2=\|r\|_2\right\}.
\ee
we have
\bo\lab{ground} If $p\in(0,2/3)$, then the following hold true
\begin{itemize}
\item [{\rm(i)}] $e\in(-\iy,0)$.
\item [{\rm(ii)}] $e$ can be achieved by $r$.
\item [{\rm(iii)}] Any minimizer of $e$ has the form as follows
$$
\left\{e^{i\theta} r(\cdot+y):\theta\in\R, y\in\R^3\right\}.
$$
\end{itemize}
\eo
\bp
{\rm(i)} Noting that $2p+2\in(2,6)$, we have
$$
\frac{1}{2p+2}=\frac{s}{2}+\frac{1-s}{6},\,\,s=\frac{3}{2p+2}-\frac{1}{2}.
$$
It follows from the interpolation inequality that there exists $C>0$ such that, for any $u\in\mathcal{M}$, $\|u\|_{2p+2}^{2p+2}\le C\|\na u\|_2^{3p}$. Then
$$
\inf_{u\in\mathcal{M}}\mathcal{E}(u)\ge\inf_{u\in\mathcal{M}}\left(\frac{1}{2}\|\na u\|_2^2-C\|\na u\|_2^{3p}\right)\ge\min_{t\in[0,\iy}\left(\frac{1}{2}t^2-Ct^{3p}\right)>-\iy.
$$
On the other hand, by the Pohozaev identity, one can get that
$$
e\le\mathcal{E}(r)=\left(\frac{3}{2}-\frac{5}{2p+2}\right)\|r\|_{2p+2}^{2p+2}-\frac{1}{4}\|\na r\|_2^4<-\frac{1}{4}\|\na r\|_2^4<0.
$$
Here we used the fact that $p\in(0,2/3)$.
\vskip0.1in
{\rm(ii)} Firstly, taking any minimization sequence $\{u_n\}$ of $e$, without loss of generality, we can assume that $u_n$ is radially symmetric and nonnegative. Since $\mathcal{E}(u_n)\rg e$ as $n\rg\iy$ and $\|u_n\|_2=\|r\|_2$, thanks to $p\in(0,2/3)$, one can show that $\{u_n\}$ is bounded in $H_{rad}^1(\R^3)$. Up to a subsequence, for some $u_0\in H_{rad}^1(\R^3)$, $u_n\rg u_0$ weakly in $H^1(\R^3)$ and strongly in $L^{2p+2}(\R^3)$ as $n\rg\iy$. If $u_0\equiv0$, then by $\mathcal{E}(u_n)\rg e$, we have
$$
\frac{1}{2}\int_{\R^3}|\na u_n|^2+\frac{1}{4}\left(\int_{\R^3}|\na u_n|^2\right)^2\rg e<0,\,\,n\rg\iy,
$$
which is a contradiction. Now, we claim that $\|u_0\|_2=\|r\|_2$. Obviously, $\|u_0\|_2\le\|r\|_2$ and $\mathcal{E}(u_0)\le e$. Then to show $e$ can be achieved by $u_0$, it suffices to rule out the case: $\|u_0\|_2<\|r\|_2$. If such case occurs, let
$$
w(\cdot)=\frac{1}{s}u_0\left(\frac{1}{t}\cdot\right),\,\,s,t\ge0,\,\,s^{2p+2}=t^3,
$$
then choosing $t>0$ such that $w\in\mathcal{M}$, i. e.,
$$
\|w\|_2^2=\frac{t^3}{s^2}\|u_0\|_2^2=\|r\|_2^2.
$$
And we have $s^{2p+2}=t^3>s^2$, which implies that $s>1$ and $s^2>t$ since $p\in(0,2/3)$. Thus, $u_0(\cdot)=sw(t\cdot)$ and
$$
\mathcal{E}(u_0)=\frac{s^2}{2t}\int_{\R^3}|\na w|^2+\frac{s^4}{4t^2}\left(\int_{\R^3}|\na w|^2\right)^2-\frac{1}{p+1}\int_{\R^3}|w|^{2p+2}>\mathcal{E}(w),
$$
which contradicts the fact that $\mathcal{E}(w)\ge e$.

Secondly, we show that $u_0=r$. Similar to \cite{Lieb}, there exists $\la_0>0$ such that
$$
-\frac{1}{2}\left(1+\int_{\R^3}|\na u_0|^2\right)\DD u_0+\la_0u_0=u_0^{2p+1},\,\,x\in\R^3.
$$
Similar as in \cite{Li}, by \cite{Kwong}, let $Q$ be the unique radially symmetric solution of
$$
-\DD Q+\la_0 Q=Q^{2p+1},\,\,Q>0,\,\,Q\in H^1(\R^3),
$$
then we have
$$
u_0(x)=Q\left(\frac{x}{\sqrt{c}}\right),\,\, \sqrt{c}=\frac{1}{2}\left(\frac{1}{2}\|\na Q\|_2^2+\sqrt{\frac{1}{4}\|\na Q\|_2^4+2}\right).
$$
Similarly,
$$
r(x)=\ti{Q}\left(\frac{x}{\sqrt{d}}\right),\,\, \sqrt{d}=\frac{1}{2}\left(\frac{1}{2}\|\na \ti{Q}\|_2^2+\sqrt{\frac{1}{4}\|\na \ti{Q}\|_2^4+2}\right),
$$
where $\ti{Q}$ is the unique radially symmetric solution of
$$
-\DD \ti{Q}+\ti{Q}=\ti{Q}^{2p+1},\,\,\ti{Q}>0,\,\,\ti{Q}\in H^1(\R^3).
$$
Let
$$
Q(\cdot)=\la_0^{\frac{1}{2p}}\bar{Q}(\la_0^{\frac{1}{2}}\cdot),
$$
then
$$
-\DD \bar{Q}+\bar{Q}=\bar{Q}^{2p+1},\,\,\bar{Q}>0,\,\,\bar{Q}\in H^1(\R^3).
$$
Then we know $\bar{Q}\equiv\ti{Q}$ and
\begin{align*}
&\|u_0\|_2^2=c^{\frac{3}{2}}\|Q\|_2^2=c^{\frac{3}{2}}\la_0^{\frac{1}{p}-\frac{3}{2}}\|\bar{Q}\|_2^2\\
&=\frac{1}{8}\left(\frac{1}{2}\la_0^{\frac{1}{p}-\frac{1}{2}}\|\na \bar{Q}\|_2^2+\sqrt{\frac{1}{4}\la_0^{\frac{2}{p}-1}\|\na \bar{Q}\|_2^4+2}\right)^3\la_0^{\frac{1}{p}-\frac{3}{2}}\|\bar{Q}\|_2^2.
\end{align*}
Since
\begin{align*}
\|r\|_2^2=\frac{1}{8}\left(\frac{1}{2}\|\na \ti{Q}\|_2^2+\sqrt{\frac{1}{4}\|\na \ti{Q}\|_2^4+2}\right)^3\|\ti{Q}\|_2^2
\end{align*}
and $\|u_0\|_2=\|r\|_2$, we get that $\la_0=1$, where we used the fact that $p\in(0,2/3)$. Therefore, $u_0$ is a radially symmetric positive solution of problem \re{lim} and we get the claim as desired.
\vskip0.1in
{\rm(iii)} The proof is similar to \cite[Theorem II.1]{Lions}. So we omit the details.
\ep

\subsection{The linearized problem}
Let $L$ be the linearization of \re{lim} at $r$ acting on $L^2(\R^3,\mathbb{C})$ with domain in $H^2(\R^3,\mathbb{C})$. Precisely, for any $\xi\in H^2(\R^3,\mathbb{C})$,
$$
L\xi=-\frac{1}{2}\left(1+\int_{\R^3}|\na r|^2\right)\DD\xi-\frac{1}{2}\left(\int_{\R^3}\na r\na(\xi+\bar{\xi})\right)\DD r+\xi-r^{2p}\left[p(\xi+\bar{\xi})+\xi\right],
$$
and
$$
L\xi=L_+\mathfrak{Re}(\xi)+iL_-\mathfrak{Im}(\xi).
$$
If $\eta\in H^2(\R^3,\R)$, then
$$
L_-\eta=-\frac{1}{2}\left(1+\int_{\R^3}|\na r|^2\right)\DD\eta+\eta-r^{2p}\eta.
$$
It is easy to check that $L_+,L_-$ are self-adjoint. Recalling that $L_-r=0$, we know $r$ is an eigenfunction of the operator
$$
-\frac{1}{2}\left(1+\int_{\R^3}|\na r|^2\right)\DD+1
$$
in $L^2(\R^3,r^{2p}\ud x)$. Since $r(x)>0,\,\, x\in\R^3$, we know $1$ is the first eigenvalue which is simple and the associated eigenfunction is $r$. Then
$$
{\rm Ker}L_-={\rm span}\{r\},
$$
and $\lan L_- \eta,\eta\ran\ge0$ for any $\eta\in H^1(\R^3)$.

\s{Proof of Theorem \ref{Th1}}

In this section, we are in position to investigate the modulational stability of ground states to problem \re{lim}. 

\subsection{Spectral estimates of $L_\pm$} To start the proof, we give some crucial lemmas as follows.

\bl{\rm\cite{Marco}}\lab{Marco1}
For any $\phi\in H^1(\R^3,\mathbb{C})$ with $\|\phi\|_2=\|r\|_2$ and
$$
\inf_{(x,\theta)\in\R^3\times[0,2\pi)}\|\phi-e^{i\theta}r(\cdot-x)\|\le\|r\|,
$$
then the minimization problem
$$
\inf_{(x,\theta)\in\R^3\times[0,2\pi)}\|\phi-e^{i\theta}r(\cdot-x)\|
$$
is achieved at some $(x_0,\gamma)\in\R^3\times[0,2\pi)$.
\el
\br\lab{inf}
For $\phi,(x_0,\gamma)$ given above, let
$$
w:=u+iv=e^{-i\gamma}\phi(\cdot+x_0)-r(\cdot),
$$
then
$$
\|w\|\le\|r\|,\,\,\|w+r\|_2=\|r\|_2.
$$
We claim that
\be\lab{ortho}
(v,r)_{H^1}=(u,\pl_{x_j}r)_{H^1}=0,\,\,j=1,2,3.
\ee
In fact, for any $(x,\theta)\in\R^3\times[0,2\pi)$, consider the function
\begin{align*}
&\Upsilon(x,\theta)=\|\phi-e^{i\theta}r(\cdot-x)\|^2\\
&=\|\phi\|^2+\|r\|^2-2\mathfrak{Re}\int_{\R^3}e^{i\theta}\bar{\phi}(y)\left(-\frac{1}{2}\DD r+r\right)(y-x)\,\ud y\\
\end{align*}
So
$$
\pa_{\theta}\Upsilon(x_0,\gamma)=\pa_{x_j}\Upsilon(x_0,\gamma)=0,\,\,j=1,2,3.
$$
Since
$$
\pa_{\theta}\Upsilon(x_0,\gamma)=-2\int_{\R^3}v\left(-\frac{1}{2}\DD r+r\right)=-2(v,r)_{H^1},
$$
and for $j=1,2,3$,
$$
\pa_{x_j}\Upsilon(x_0,\gamma)=-2\int_{\R^3}u\left[-\frac{1}{2}\DD(\pa_{x_j}r)+\pa_{x_j}r\right]=-2(u,\pl_{x_j}r)_{H^1},
$$
Thus, we get \re{ortho}.
\er
Similar to \cite[Lemma 2.1]{Marco}, let $$\mathcal{V}:=\left\{u\in H^1(\R^3):(u,r)=0\right\},$$ then we have
\bl
$
\inf_{u\in\mathcal{V}}\lan L_+u,u\ran=0.
$
\el
Set
$$
\mathcal{V}_0:=\left\{u\in H^1(\R^3):(u,r)=(u,\pa_{x_j}r)_{H^1}=0,\,j=1,2,3\right\},
$$
we have
\bl
$
\inf_{u\in\mathcal{V}_0}\frac{\lan L_+u,u\ran}{\|u\|^2}>0.
$
\el
\bp
It suffices to show that
\be\lab{positive}
\inf_{u\in\mathcal{V}_0}\frac{\lan L_+u,u\ran}{\|u\|_2^2}>0.
\ee
Indeed, if \re{positive} holds true, then we have
$$
\inf_{u\in\mathcal{V}_0}\frac{\lan L_+u,u\ran}{\|u\|^2}>0.
$$
If not, there exists $\{u_n\}\subset\mathcal{V}_0$ satisfying $\|u_n\|=1$ and $\lan L_+u_n,u_n\ran\rg0$ as $n\rg\iy$. By \re{positive}, $u_n\rg0$ strongly in $L^2(\R^3)$ and then weakly in $H^1(\R^3)$ as $n\rg\iy$. So
\begin{align*}
\lan L_+u_n,u_n\ran&=\frac{1}{2}\left(1+\int_{\R^3}|\na r|^2\right)\|\na u_n\|_2^2+\left(\int_{\R^3}\na r\na u_n\right)^2+\int_{\R^3}\left[1-(2p+1)r^{2p}\right]u_n^2\\
&=\frac{1}{2}\left(1+\int_{\R^3}|\na r|^2\right)\|\na u_n\|_2^2+o_n(1),\,\,\mbox{as}\,\,n\rg\iy.
\end{align*}
It yields that $u_n\rg0$ strongly in $H^1(\R^3)$ as $n\rg\iy$, which contradicts the fact that $\|u_n\|=1$ for any $n$.

In the following, we only need to show \re{positive} is true. If not, there exists $\{u_n\}\subset \mathcal{V}_0$ with $\|u_n\|_2=1$ such that $\lan L_+u_n,u_n\ran\rg0$ as $n\rg\iy$. Noting that
$$
\lan L_+u_n,u_n\ran=\frac{1}{2}\left(1+\int_{\R^3}|\na r|^2\right)\|\na u_n\|_2^2+\left(\int_{\R^3}\na r\na u_n\right)^2+\|u_n\|_2^2-(2p+1)\int_{\R^3}r^{2p}u_n^2,
$$
we get
\begin{align}\lab{cov1}
&\limsup_{n\rg\iy}\left[\frac{1}{2}\left(1+\int_{\R^3}|\na r|^2\right)\|\na u_n\|_2^2+\|u_n\|_2^2\right]\nonumber\\
&=\limsup_{n\rg\iy}\left[(2p+1)\int_{\R^3}r^{2p}u_n^2-\left(\int_{\R^3}\na r\na u_n\right)^2\right],
\end{align}
and
$
\limsup_{n\rg\iy}\|\na u_n\|_2^2\le2(2p+1)\|r\|_\iy^{2p}.
$
So, $\{u_n\}$ is bounded in $H^1(\R^3)$. Up to a subsequence, there exists $u\in H^1(\R^3)$ such that $u_n\rg u$ weakly in $H^1(\R^3)$ and a. e. in $\R^3$ as $n\rg\iy$. Obviously, $u\in\mathcal{V}_0$ and $\lan L_+u,u\ran\ge0$. On the other hand, since $r(x)\rg0$ as $|x|\rg\iy$, up to a subsequence, $\int_{\R^3}r^{2p}u_n^2\rg\int_{\R^3}r^{2p}u^2$ as $n\rg\iy$. Thanks to the weak lower semi-continuity of the norm, we have
$$
\lan L_+u,u\ran\le\liminf_{n\rg\iy}\lan L_+u_n,u_n\ran=0,
$$
which implies that $\lan L_+u,u\ran=0$. That is,
\begin{align}\lab{cov2}
\frac{1}{2}\left(1+\int_{\R^3}|\na r|^2\right)\|\na u\|_2^2+\|u\|_2^2=(2p+1)\int_{\R^3}r^{2p}u^2-\left(\int_{\R^3}\na r\na u\right)^2.
\end{align}
Noting that
$$
\lim_{n\rg\iy}\int_{\R^3}\na r\na u_n=\int_{\R^3}\na r\na u,
$$
it follows from \re{cov1} and \re{cov2} that $u_n\rg u$ strongly in $H^1(\R^3)$ as $n\rg\iy$ and $\|u\|_2=1$. Then there exist Lagrange multipliers $\la,\mu,\la_1,\la_2,\la_3$ such that for any $\eta\in H^1(\R^3)$,
\be\lab{lagr}
\lan L_+u,\eta\ran=\la(u,\eta)+\mu(r,\eta)+\sum_{i=1}^3\la_i(\pl_{x_i}r,\eta)_{H^1}.
\ee
Thanks to $u\in\mathcal{V}_0$, $\la=0$. For any $j$,
$$
\lan L_+u,\pl_{x_j}r\ran=\lan L_+\pl_{x_j}r,u\ran=0,
$$
where we used the fact that
$
{\rm Ker}L_+={\rm span}\{\pl_{x_1}r,\pl_{x_2}r,\pl_{x_3}r\}.
$
Then taking $\eta=\pl_{x_j}r$, for $j=1,2,3$, we have
$$
\la_j(\pl_{x_j}r,\pl_{x_j}r)_{H^1}=0,
$$
and $\la_j=0$. Here we used the fact that
$
(r,\pl_{x_j}r)=0,\,\,j=1,2,3,
$
and
$
(\pl_{x_i}r,\pl_{x_j}r)_{H^1}=0,\,\,i\not=j.
$
In turn, for any $\eta\in H^1(\R^3)$, it holds true that $\lan L_+u,\eta\ran=\mu(r,\eta)$. If $\mu=0$, then $u\in{\rm Ker}L_+$, which contradicts the fact that $\|u\|_2=1$ and $(u,\pa_{x_j}r)_{H^1}=0,\,j=1,2,3$. So $\mu\not=0$.

In the following, we show that we can reach a contradiction: $\mu=0$. In fact, by computation, one can get that
$$
\DD(x\cdot\na r)=2\DD r+\sum_{j=1}^3x_j\pl_{x_j}\DD r,\,\,\int_{\R^3}\na r\na(x\cdot\na r)=-\frac{1}{2}\|\na r\|_2^2.
$$
Then
\begin{align}\lab{zuhe1}
L_+(x\cdot\na r)&=-\left(1+\frac{1}{2}\|\na r\|_2^2\right)\DD r-\sum_{j=1}^3x_j\left[-\frac{1}{2}(1+\|\na r\|_2^2)\DD r+r-r^{2p+1}\right]\nonumber\\
&=-\left(1+\frac{1}{2}\|\na r\|_2^2\right)\DD r.
\end{align}
Meanwhile, since $r$ is a solution of \re{lim}, we have
\begin{align}\lab{zuhe2}
L_+r&=-\frac{1}{2}(1+3\|\na r\|_2^2)\DD r+r-(2p+1)r^{2p+1}\nonumber\\
&=\left[p+(p-1)\|\na r\|_2^2\right]\DD r-2pr.
\end{align}
So by \re{zuhe1}-\re{zuhe2}, we get that
$$
L_+\left(\frac{r}{2p}+\frac{p+(p-1)\|\na r\|_2^2}{p(2+\|\na r\|_2^2)}(x\cdot\na r)\right)=-r.
$$
Recalling that $L_+u=\mu r$, for some $\vartheta\in\R^3$, we have
$$
u=\vartheta\cdot\na r-\mu\left[\frac{r}{2p}+\frac{p+(p-1)\|\na r\|_2^2}{p(2+\|\na r\|_2^2)}(x\cdot\na r)\right].
$$
Thanks to the fact that
$$
\int_{\R^3}ru=\int_{\R^3}r\pl_{x_j}r=0,\,\,j=1,2,3,
$$
we reach that
$$
\mu\int_{\R^3}\left[\frac{r^2}{2p}+\frac{p+(p-1)\|\na r\|_2^2}{p(2+\|\na r\|_2^2)}(x\cdot\na r)r\right]=0.
$$
Since
$$
\int_{\R^3}(x\cdot\na r)r=-\frac{3}{2}\int_{\R^3}r^2,
$$
then by $p\in(0,2/3)$,
\begin{align*}
\int_{\R^3}\left[\frac{r^2}{2p}+\frac{p+(p-1)\|\na r\|_2^2}{p(2+\|\na r\|_2^2)}(x\cdot\na r)r\right]=\frac{(2-3p)+(4-3p))\|\na r\|_2^2}{2p(2+\|\na r\|_2^2)}\int_{\R^3}r^2>0.
\end{align*}
It follows that $\mu=0$. The proof is complete.
\ep

\bl\lab{zui}
Let $w=u+iv\in H^1(\R^3,\mathbb{C})$ with $u,v\in H^1(\R^3)$. If $\|w+r\|_2=\|r\|_2$ and
$$
(u,\pa_{x_j}r)_{H^1}=0,\,j=1,2,3,
$$
then there exist $D,D_1,D_2$ such that
$$
\lan L_+u,u\ran\ge D\|u\|^2-D_1\|w\|^4-D_2\|w\|^3.
$$
\el
\bp
By $\|w+r\|_2=\|r\|_2$, we get $(u,r)=-\frac{1}{2}\|w\|_2^2$. Without loss of generality, we assume that $\|r\|_2=1$. Let
$$
u=u_{\parallel}+u_{\perp},\,\,u_{\parallel}=(u,r)r,
$$
then $(u_{\perp},r)=0$. Noting that $(r,\pa_{x_j}r)_{H^1}=0,\,j=1,2,3$, we have $(u_{\perp},\pa_{x_j}r)_{H^1}=0,\,j=1,2,3$ and $u_{\perp}\in\mathcal{V}_0$. It follows that $\lan L_+u_{\perp},u_{\perp}\ran\ge C\|u_{\perp}\|_2^2$ for some $C>0$. Similar to \cite{Marco},
\be\lab{low1}
\lan L_+u_{\perp},u_{\perp}\ran\ge C(\|u\|^2-\|w\|_2^4).
\ee

On the other hand, since $r$ is a solution of \re{lim}, we have
$$
\frac{1}{2}(1+\|\na r\|_2^2)\|\na r\|_2^2+\|r\|_2^2=\|r\|_{2p+2}^{2p+2}.
$$
It follows that
$$
\lan L_+r,r\ran=-p\|\na r\|_2^2-2p\|r\|_2^2+(1-p)\|\na r\|_2^4\ge-2p\|r\|^2,
$$
and then
\be\lab{low2}
\lan L_+u_{\parallel},u_{\parallel}\ran=\frac{1}{4}\|w\|_2^4\lan L_+r,r\ran\ge-\frac{p}{2}\|w\|_2^4\|r\|^2.
\ee
Finally, since $r$ satisfies \re{lim}, we get
$$
L_+r=p(1+\|\na r\|_2^2)\DD r-\|\na r\|_2^2\DD r-2pr.
$$
Then for some $C>0$,
\begin{align*}
\lan L_+u_{\perp},r\ran&=\lan L_+r,u_{\perp}\ran=\left[(1-p)\|\na r\|_2^2-p\right]\int_{\R^3}\na r\na u_{\perp}\\
&=\left[(1-p)\|\na r\|_2^2-p\right]\left(\int_{\R^3}\na r\na u-\int_{\R^3}\na r\na u_{\parallel}\right)\\
&=\left[(1-p)\|\na r\|_2^2-p\right]\left(\int_{\R^3}\na r\na u+\frac{1}{2}\|w\|_2^2\int_{\R^3}|\na r|^2\right)\\
&\le C(\|\na u\|_2+\|w\|_2^2).
\end{align*}
So
\begin{align}\lab{low3}
\lan L_+u_{\perp},u_{\parallel}\ran&=-\frac{1}{2}\|w\|_2^2\lan L_+u_{\perp},r\ran\ge-\frac{C}{2}\|w\|_2^2(\|\na u\|_2+\|w\|_2^2)\nonumber\\
&\ge-C(\|w\|^3+\|w\|^4).
\end{align}
Thus, the result as claimed is yielded by \re{low1}-\re{low3}.
\ep

\bl\lab{minus}
$$
\inf_{\stackrel{v\in H^1(\R^3)\setminus\{0\}}{(v,r)_{H^1}=0}}\frac{\lan L_-v,v\ran}{\|v\|^2}>0.
$$
\el
\bp
It suffices to show that
$$
\omega:=\inf_{\stackrel{v\in H^1(\R^3)\setminus\{0\}}{\|v\|_2=1,\,(v,r)_{H^1}=0}}\lan L_-v,v\ran>0.
$$
Since $r(x)\rg0$ as $|x|\rg\iy$, similar to \cite{Marco}, we know $\omega\ge0$. If $\omega=0$, taking any minimizing sequence $\{v_n\}$, $\{v_n\}$ is bounded in $H^1(\R^3)$ and for some $v\in H^1(\R^3)$, we have $v_n\rg v$ weakly in $H^1(\R^3)$ and a. e. in $\R^3$ as $n\rg\iy$. So $(v,r)_{H^1}=0$ and by the decay of $r$, $$
\lim_{n\rg\iy}\int_{\R^3}r^{2p}v_n^2=\int_{\R^3}r^{2p}v^2.
$$
$$
0\le\lan L_-v,v\ran\le\liminf_{n\rg\iy}\lan L_-v_n,v_n\ran=0
$$
Then $\lan L_-v,v\ran=0$. Furthermore, we know $v_n\rg v$ strongly in $H^1(\R^3)$ as $n\rg\iy$ and $\|v\|_2=1$. In turn, there exist $\la,\mu$ such that
$$
\lan L_-v,\eta\ran=\la(v,\eta)+\mu(r,\eta)_{H^1},\,\,\eta\in H^1(\R^3).
$$
By taking $\eta=v$, $\la=0$. Finally, we take $\eta=r$ and get that
$$
\mu\|r\|_{H^1}^2=\lan L_-v,r\ran=\lan L_-r,v\ran=0.
$$
That is, $\mu=0$ and $L_-v=0$. Recalling that $\mbox{Ker}L_-={\rm Span}\{r\}$, we get that $v=\theta r$ for some $\theta\in\R$. Noting that $(v,r)_{H^1}=0$, $\theta=0$, which contradicts the fact that $\|v\|_2=1$.
\ep

\subsection{Proof of Theorem \ref{Th1}}
\bp
Take $\phi=r+w,(x_0,\gamma)$ given in Lemma \ref{Marco1} and $w,u,v$ given in Remark \ref{inf}. Let $I(\phi)=\mathcal{E}(\phi)+\|\phi\|_2^2$, we get that
$I'(r)=0$ in $H^{-1}(\R^3)$ and then by Proposition \ref{ground}, $I(\phi)\ge I(r)$. By the Taylor expend, for some $\theta\in(0,1)$, we have
\begin{align*}
I(\phi)-I(r)&=I(r+w)-I(r)=\frac{1}{2}\lan I''(r+\theta w)w,w\ran\\
&:=\lan L_+u,u\ran+\lan L_-v,v\ran+J+K,
\end{align*}
where
\begin{align*}
J=&\frac{1}{2}\left(\|\na(r+\theta u)\|_2^2+\theta^2\|\na v\|_2^2-\|\na r\|_2^2\right)\left(\|\na u\|_2^2+\|\na v\|_2^2\right)\\
&+\left[\theta\|\na u\|_2^2+\theta\|\na v\|_2^2+\int_{\R^3}\na r\na u\right]^2-\left(\int_{\R^3}\na r\na u\right)^2,
\end{align*}
and
\begin{align*}
K=&\int_{\R^3}[(2p+1)r^{2p}u^2+r^{2p}v^2-|w|^2|r+\theta w|^{2p}]\\
&-2p\int_{\R^3}|r+\theta w|^{2p-2}|ru+\theta u^2+\theta v^2|^2.
\end{align*}
It is easy to check that $J\ge-C(\|w\|^3+\|w\|^4)$ for some $C>0$. Similar as that in \cite{Wein2}, by an interpolation estimate of Nirenberg and Gagliardo, one can get that $K\ge-C(\|w\|^{2+\tau}+\|w\|^6)$, where $\tau>0$ and $C>0$. Finally, the claim is concluded by \re{ortho}, Remark \ref{inf}, Lemma \ref{zui} and Lemma \ref{minus}.
\ep

\end{document}